\documentclass{amsart}
\usepackage{amssymb,amscd, supertabular}

\DeclareMathSymbol{\twoheadrightarrow}  {\mathrel}{AMSa}{"10}

\def\GG{{\mathcal G}}
\def\HH{{\mathcal H}}

\def\Q{{\mathbb Q}}
\def\Z{{\mathbb Z}}
\def\C{{\mathbb C}}

\def\F{{\mathbb F}}
\def\P{{\mathbb P}}

\def\RR{{\mathfrak R}}

\def\Perm{\mathrm{Perm}}
\def\Gal{\mathrm{Gal}}
\def\PSL{\mathrm{PSL}}
\def\PGL{\mathrm{PGL}}
\def\End{\mathrm{End}}
\def\Aut{\mathrm{Aut}}
\def\Hom{\mathrm{Hom}}

\def\SS{{\mathbf S}}

\def\fchar{\mathrm{char}}
\def\GL{\mathrm{GL}}

\def\SL{\mathrm{SL}}
\def\Sp{\mathrm{Sp}}

\def\M{\mathrm{M}}

\def\dim{\mathrm{dim}}

\def\O{{\mathcal O}}

\def\a{{\mathfrak a}}
\def\b{{\mathfrak b}}

           \def\hh{{\mathbf h}}

\DeclareMathOperator{\mG}{G}
\DeclareMathOperator{\mL}{L}

\newtheorem{thm}{Theorem}[section]
\newtheorem{lem}[thm]{Lemma}
\newtheorem{cor}[thm]{Corollary}
\theoremstyle{definition}
\newtheorem{defn}[thm]{Definition}
\newtheorem{ex}[thm]{Example}
\newtheorem{rem}[thm]{Remark}

\title[Hyperelliptic jacobians]
{Endomorphism algebras of hyperelliptic jacobians and finite
projective lines}

\author[Arsen Elkin]{Arsen Elkin}
\address{Institute of Mathematics, Hebrew University of Jerusalem, Givat Ram, Jerusalem, 91904, Israel} \email{arsen\char`\@math.huji.ac.il}
\author[Yuri\ G.\ Zarhin]{Yuri\ G.\ Zarhin}
\address{Department of Mathematics, Pennsylvania State University,
University Park, PA 16802, USA} \email{zarhin\char`\@math.psu.edu}

\begin{document}
\maketitle

\section{Statement of results}

 Let $K$ be a field with $\fchar(K)\ne 2$. Let us fix an algebraic closure
$K_a$ of $K$. Let us put $\Gal(K):=\Aut(K_a/K)$. If $X$ is an
abelian variety of positive dimension
 over $K_a$ then we write $\End(X)$ for the ring of all its
$K_a$-endomorphisms and $\End^0(X)$ for the corresponding
(semisimple finite-dimensional)  $\Q$-algebra $\End(X)\otimes\Q$.
We write $\End_K(X)$ for the ring of all  $K$-endomorphisms of $X$
and $\End_K^0(X)$ for the corresponding (semisimple
finite-dimensional) $\Q$-algebra $\End_K(X)\otimes\Q$. The
absolute Galois group $\Gal(K)$ of $K$ acts on $\End(X)$ (and
therefore on $\End^0(X)$) by ring (resp. algebra) automorphisms and
$$\End_K(X)=\End(X)^{\Gal(K)}, \ \End_K^0(X)=\End^0(X)^{\Gal(K)},$$
since every endomorphism of $X$ is defined over a finite separable
extension of $K$.

If $n$ is a positive integer that is not divisible by $\fchar(K)$
then we write $X_n$ for the kernel of multiplication by $n$ in
$X(K_a)$. It is well-known \cite{MumfordAV} that  $X_n$ is a free
$\Z/n\Z$-module of rank $2\dim(X)$. In particular, if $n=\ell$ is
a prime then $X_{\ell}$ is an $\F_{\ell}$-vector space of
dimension $2\dim(X)$.

 If $X$ is defined over $K$ then $X_n$ is a Galois
submodule in $X(K_a)$. It is known that all points of $X_n$ are
defined over a finite separable extension of $K$. We write
$\bar{\rho}_{n,X,K}:\Gal(K)\to \Aut_{\Z/n\Z}(X_n)$ for the
corresponding homomorphism defining the structure of the Galois
module on $X_n$,
$$\tilde{G}_{n,X,K}\subset
\Aut_{\Z/n\Z}(X_{n})$$ for its image $\bar{\rho}_{n,X,K}(\Gal(K))$
and $K(X_n)$ for the field of definition of all points of $X_n$.
Clearly, $K(X_n)$ is a finite Galois extension of $K$ with Galois
group $\Gal(K(X_n)/K)=\tilde{G}_{n,X,K}$. 
If $n=\ell$  then we get a natural faithful linear representation
$$\tilde{G}_{\ell,X,K}\subset \Aut_{\F_{\ell}}(X_{\ell})$$
of $\tilde{G}_{\ell,X,K}$ in the $\F_{\ell}$-vector space
$X_{\ell}$. Recall \cite{Silverberg} that all endomorphisms of $X$
are defined over $K(X_4)$; this gives rise to the natural
homomorphism $$\kappa_{X,4}:\tilde{G}_{4,X,K,} \to
\Aut(\End^0(X))$$ and $\End_K^0(X)$ coincides with the subalgebra
$\End^0(X)^{\tilde{G}_{4,X,K}}$ of $\tilde{G}_{4,X,K}$-invariants
\cite[Sect. 1]{ZarhinLuminy}.

Let $f(x)\in K[x]$ be a polynomial of degree $n\ge 3$ without
multiple roots. Let $\RR_f\subset K_a$ be the $n$-element set of
roots of $f$. Then $K(\RR_f)$ is the splitting field of $f$ and
$\Gal(f):=\Aut(K(\RR_f)/K)$ is the Galois group of $f$ (over $K$).
One may view $\Gal(f)$ as a group of permutations of $\RR_f$; it
is transitive if and only if $f(x)$ is irreducible.

Let us  consider   the hyperelliptic curve $C_f:y^2=f(x)$ and its
jacobian $J(C_f)$. It is well-known \cite{ZarhinTexel} that $J(C_f)$ is a
$\left[\frac{n-1}{2}\right]$-dimensional abelian variety defined
over $K$. The aim of this paper is to study $\End^0(J(C_f))$,
assuming that $n=q+1$ where $q$ is a  power of a prime $p$ and
$\Gal(f)=\PSL_2(\F_q)$ acts via fractional-linear transformations
on $\RR_f$ identified with the projective line $\P^1(\F_q)$. It
follows from results of \cite{ZarhinMRL,ZarhinMMJ,ZarhinBSMF} that
for every $q$ in all characteristics there
exist $K$ and $f$ with
$\End^0(J(C_f))=\Q$. On the other hand, it is known that
$\End^0(J(C_f))=\Q$ \cite{ZarhinMRL,ZarhinTexel, ZarhinSbg} if
$p=2$ and $q \ge 8$. (It is also true for $q=4$ if one assumes
that $\fchar(K)\ne 3$ \cite{ZarhinBSMF}.) However, if $q=5$ then
there are examples where $\End^0(J(C_f))$ is a (real) quadratic
field (even in characteristic zero)
\cite{BruRankJ,HashiOnBru,WilEM,Elkin}.

Our main result is the following statement.

\begin{thm}
\label{main} Let us assume that $\fchar(K)=0$.
 Suppose that $n=q+1$ where $q\ge 5$ is a prime  power that is
congruent to $\pm 3$ modulo $8$. Suppose that $f(x)$ is
irreducible and $\Gal(f)\cong \PSL_2(\F_q)$. Then one of the
following two conditions holds:
\begin{itemize}
 \item[(i)]
$\End^0(J(C_f))=\Q$ or a quadratic field. In particular, $J(C_f)$
is an absolutely simple abelian variety.

\item [(ii)] $q$ is congruent to $3$ modulo $8$ and $J(C_f)$ is
$K_a$-isogenous to a self-product of an elliptic curve with
complex multiplication by $\Q(\sqrt{-q})$.
\end{itemize}
\end{thm}

\begin{rem}
It follows from results of \cite{Katsura} (see also \cite{Shioda},
\cite{Schoen})  that if the case (ii) of Theorem \ref{main} holds
then $J(C_f)$ is isomorphic over $K_a$ to a product of mutually
isogenous elliptic curves with complex multiplication by
$\Q(\sqrt{-q})$.
\end{rem}

The paper is organized as follows.  Section \ref{luminy} contains
auxiliary results about endomorphism algebras of abelian
varieties. In Section \ref{mainp} we prove the main result. In
Section \ref{noniso} (and Section \ref{quadratic}) we prove the
absolute simplicity of $J(C_f)$ when $q\ge 11$ is congruent to $3$
modulo $8$ and $K=\Q$.  Section \ref{examples} contains examples.

\section{Endomorphism algebras of abelian varieties}
\label{luminy}

\begin{rem}
\label{minimal}
Recall \cite{FT} (see also \cite[p. 199]{ZarhinP})
that a surjective homomorphism of finite groups
$\pi:\GG_1\twoheadrightarrow \GG$ is called a {\sl minimal cover}
if no proper subgroup of $\GG_1$ maps onto $\GG$. If $H$ is a
normal subgroup of $\GG_1$ that lies in $\ker(\pi)$  then the
induced surjection $\GG_1/H\twoheadrightarrow \GG$ is also a
minimal cover.

\begin{itemize}
\item[(i)] If a surjection $\GG_2\twoheadrightarrow \GG_1$ is also
a minimal cover then one may easily check that the composition
$\GG_2 \to\GG$ is surjective and a minimal cover.

\item[(ii)] Clearly, if $\GG$ is simple then every proper normal
subgroup in $\GG_1$  lies in $\ker(\pi)$.

\item[(iii)] If $\GG$ is perfect then its universal central
extension is a minimal cover \cite{Suzuki}.

\item[(iv)] If $\GG'\twoheadrightarrow \GG$ is an arbitrary
surjective homomorphism of finite groups then
 there always exists a subgroup $\HH\subset \GG'$ such that
$\HH\to\GG$ is surjective and a minimal cover. Clearly, if $\GG$
is perfect then $\HH$ is also perfect.

\end{itemize}
\end{rem}

 The field inclusion
$K(X_2)\subset K(X_4)$ induces a natural surjection \cite[Sect.
1]{ZarhinLuminy}
$$\tau_{2,X}:\tilde{G}_{4,X,K}\to\tilde{G}_{2,X,K}.$$

\begin{defn} We say that $K$ is 2-{\sl balanced} with respect to
$X$ if $\tau_{2,X}$ is a minimal cover.
\end{defn}

\begin{rem}
\label{overL}
 Clearly, there always exists a subgroup $H
\subset \tilde{G}_{4,X,K}$ such that $H\to\tilde{G}_{2,X,K}$ is
surjective and a minimal cover. Let us put $L=K(X_4)^H$. Clearly,
$$K \subset L \subset K(X_4), \ L\bigcap K(X_2)=K$$
and $L$ is a maximal overfield of $K$ that enjoys these
properties. It is also clear that
$$K(X_2)\subset L(X_2),\ L(X_4)=K(X_4), \ H=\tilde{G}_{4,X,L},
\tilde{G}_{2,X,L}=\tilde{G}_{2,X,K}$$
 and $L$ is $2$-{\sl balanced}   with respect to
$X$.
\end{rem}

The following assertion (and its proof) is (are) inspired by
Theorem 1.6 of \cite{ZarhinLuminy} (and its proof).

\begin{thm}
\label{rep}
 Suppose that $E:=\End_K^0(X)$ is a field that contains the
center $C$ of $\End^0(X)$. Let $C_{X,K}$ be the centralizer of
$\End_K^0(X)$ in $\End^0(X)$.

Then:

\begin{itemize}
\item[(i)] $C_{X,K}$ is a central simple $E$-subalgebra in
$\End^0(X)$. In addition, the centralizer of $C_{X,K}$ in
$\End^0(X)$ coincides with $E=\End_K^0(X)$ and
$$\dim_E(C_{X,K})=\frac{\dim_C(\End^0(X))}{[E:C]^2}.$$
 \item[(ii)] Assume that $K$ is $2$-balanced with
respect to $X$ and $\tilde{G}_{2,X,K}$ is a non-abelian simple
group. If $\End^0(X)\ne E$ (i.e., not all endomorphisms of $X$ are
defined over $K$) then there exist a finite perfect group $\Pi \subset
C_{X,K}^{*}$ and a surjective homomorphism $\Pi \to
\tilde{G}_{2,X,K}$ that is a minimal cover. In addition, the
induced homomorphism
$$E[\Pi] \to C_{X,K}$$
is surjective, i.e., $C_{X,K}$ is isomorphic to a direct summand of the group
algebra $E[\Pi]$.
\end{itemize}
\end{thm}

\begin{proof}
Since $E$ is a field, $C$ is  a subfield of $E$ and therefore
$\End^0(X)$ is a central simple $C$-algebra. Now the assertion (i)
follows from Theorem 2 of Sect. 10.2 in Chapter VIII of
\cite{Bourbaki}.

Now let us prove the assertion (ii).

Recall that there is the homomorphism
$$\kappa_{X,4}:\tilde{G}_{4,X,K} \to \Aut(\End^0(X))$$ such that
$$\End^0(X)^{\tilde{G}_{4,X,K}}=\End_K^0(X)=E\supset C.$$
This implies that $$\kappa_{X,4}(\tilde{G}_{4,X,K})\subset
\Aut_E(\End^0(X))\subset \Aut_C(\End^0(X))$$ and we get a
homomorphism
$$\kappa_E:\tilde{G}_{4,X,K} \to \Aut_E(C_{X,K})$$ such that
\begin{equation}
C_{X,K}^{\tilde{G}_{4,X,K}}=E \label{eqn1}
\end{equation}
Assume that $E = C_{X,K}$, i.e., $E$ coincides with its own
centralizer in $\End^0(X)$. It follows from the Skolem-Noether
theorem
that $\Aut_C(\End^0(X))=\End^0(X)^{*}/C^*$. This implies that the
group $$\Aut_E(\End^0(X))=C_{X,K}^*/C^*=E^*/C^*$$ is commutative.
It follows that $\kappa_{X,4}(\tilde{G}_{4,X,K})$ is commutative.
Since $\tilde{G}_{4,X,K}$ is perfect,
$\kappa_{X,4}(\tilde{G}_{4,X,K})$ is perfect commutative and
therefore trivial, i.e., $\End^0(X)=\End_K^{0}(X)$.

Assume that $E \ne C_{X,K}$. This means that the group
$\Gamma:=\kappa_E(\tilde{G}_{4,X,K})$ is not $\{1\}$, i.e.,
$\ker(\kappa_E)\ne \tilde{G}_{4,X,K}$. Clearly, $\Gamma$ is a
finite perfect subgroup of $\Aut_E(C_{X,K})$.

The minimality of $\tau_{2,X}$ and the simplicity of
$\tilde{G}_{2,X,K}$ imply the existence of a minimal cover
$$\Gamma \to \tilde{G}_{2,X,K},$$
 thanks to Remark \ref{minimal}.

Since $C_{X,K}$ is a central simple $E$-algebra, all its
automorphisms are inner, i.e., $\Aut_E(C_{X,K})=C_{X,K}^*/E^{*}$.
Let $\Delta \twoheadrightarrow \Gamma$ be the universal central
extension of $\Gamma$. It is well-known \cite[Ch. 2, Sect. 9]{Suzuki} that $\Delta$ is a finite
perfect group. The universality property implies that $\Delta
\twoheadrightarrow \Gamma$ is a minimal cover and the inclusion
map $\Gamma\subset C_{X,K}^*/E^{*}$ lifts (uniquely) to a
homomorphism $\pi:\Delta \to C_{X,K}^{*}$. Clearly, $\ker(\pi)$
lies in the kernel of $\Delta \twoheadrightarrow \Gamma$ and we
get a minimal cover
$$\pi(\Delta)\cong \Delta/\ker(\pi) \twoheadrightarrow \Gamma ,$$
thanks to Remark \ref{minimal}.
 Taking the compositions of minimal
covers $\pi(\Delta)\twoheadrightarrow \Gamma$ and $\Gamma
\twoheadrightarrow \tilde{G}_{2,X,K}$, we obtain a minimal cover
$\pi(\Delta)\twoheadrightarrow \tilde{G}_{2,X,K}$. If we put
$$\Pi:=\pi(\Delta)\subset C_{X,K}^{*}.$$
then we get  a minimal cover
$$\Pi\twoheadrightarrow\tilde{G}_{2,X,K}.$$
The equality \eqref{eqn1} means that
the centralizer of $\pi(\Delta)=\Pi$ in $C_{X,K}$ coincides with
$E$.
 It  follows  that if $E[\Pi]$ is the group $E$-algebra 
of $\Pi$ then the inclusion  $\Pi\subset C_{X,K}^{*}$ induces the
$E$-algebra homomorphism $\omega: E[\Pi] \to C_{X,K}$ such that  
the centralizer of its image  in $C_{X,K}$ coincides with $E$.

We claim that $\omega(E[\Pi])=C_{X,K}$
 and therefore $C_{X,K}$ is
isomorphic to a direct summand of $E[\Pi]$. This claim follows 
easily from the next lemma that was  proven  in \cite[Lemma
1.7]{ZarhinLuminy}

\begin{lem}
\label{sur} Let $F$ be a field of characteristic zero, $T$ a
semisimple finite-dimensional $F$-algebra, $S$ a
finite-dimensional central simple $F$-algebra, $\beta: T \to S$ an
$F$-algebra homomorphism that sends $1$ to $1$. Suppose that the
centralizer of the image $\beta(T)$ in $S$ coincides with the
center $F$. Then $\beta$ is surjective, i. e. $\beta(T)=S$.
\end{lem}
\end{proof}

\begin{thm}
\label{center} Suppose that $\End^0(X)$ is a simple $\Q$-algebra,
$\tilde{G}_{2,X,K}$ is a simple non-abelian group, whose order is
not a divisor of $2\dim(X)$ and
$\End_{\tilde{G}_{2,X,K}}(X_2)\cong \F_4$.

Then the center $C$ of  $\End^0(X)$ is either $\Q$ or a quadratic field. In
addition, there exists a finite separable field extension $L/K$ such that
$\tilde{G}_{2,X,L}=\tilde{G}_{2,X,K}$, the map $\tau_{2,X}:\tilde{G}_{4,X,L}
\twoheadrightarrow \tilde{G}_{2,X,L}$ is surjective and a minimal cover, the
$\Q$-algebra $E:=\End_L^0(X)$ is either $\Q$ or a quadratic field and $C\subset
\End_L^0(X)$. In particular, $C$ is either $\Q$ or a quadratic field. If
$\End^0(X)\ne \End_L^0(X)$ and $C_{X,L}$ is the centralizer of $E$ in
$\End^0(X)$ then there exist a finite perfect group $\Pi \subset C_{X,L}^{*}$
and a surjective homomorphism $\Pi \to \tilde{G}_{2,X,K}$ that is a minimal
cover. In addition, the induced homomorphism
$$E[\Pi] \to C_{X,L}$$
is surjective, i.e., $C_{X,L}$ is isomorphic to a direct summand of the group
algebra $E[\Pi]$.
\end{thm}

\begin{proof}
Choose a field $L$ as in Remark \ref{overL}. Then
$\tilde{G}_{2,X,L}=\tilde{G}_{2,X,K}$, the map
$$\tau_{2,X}:\tilde{G}_{4,X,L}\to
\tilde{G}_{2,X,L}=\tilde{G}_{2,X,K}$$ is surjective and a minimal
cover. We have
$$\End_L(X)\otimes \Z/2\Z \hookrightarrow
\End_{\Gal(L)}(X_2)=\End_{\tilde{G}_{2,X,L}}(X_2)=\F_4.$$ It
follows that the rank of the free $\Z$-module $\End_L(X)$ is $1$
or $2$; Lemma 1.3 of \cite{ZarhinLuminy} implies that $\End_L(X)$
has no zero divisors. This implies that
$\End_L^0(X)=\End_L(X)\otimes\Q$ is a division algebra of
$\Q$-dimension $1$ or $2$. This means that $E=\End_L^0(X)$ is
either $\Q$ or a quadratic field.

 Recall that the center $C$ of $\End^0(X)$ is
a number field, whose degree $[C:\Q]$ divides $2\dim(X)$. The
group $\tilde{G}_{4,X,L}$ acts via automorphisms on $C$ and
$$C^{\tilde{G}_{4,X,L}}=C\bigcap
\End_L^0(X)$$ is either $\Q$ or a quadratic field. Since
$\tilde{G}_{2,X,L}=\tilde{G}_{2,X,K}$ has no normal subgroups of
index dividing $2\dim(X)$ , the same is true for
$\tilde{G}_{4,X,L}$  and therefore $\tilde{G}_{4,X,L}$ acts on $C$
trivially, i.e., $C \subset \End_L^0(X)$.
 In order to finish the proof, one has only to apply
Theorem \ref{rep} ( to $L$ instead of $K$).
\end{proof}


\begin{defn}
We say that a group is {\em FTKL-exceptional} if it is one of the following:
\begin{itemize}
\item[(i)] $\Sp_{2n}(q)$ for some even $q$ and $n\geq 2$, except $\Sp_4(2)'$ and $\Sp_6(2)$;
\item[(ii)] $\Omega_{2n}^\pm(\F_q)$ for some even $q$ and $n\geq 4$, except $\Omega_8^+(2)$;
\item[(iii)] $\mL_4(q)$ for some even $q$, except $\mL_4(2)$; or
\item[(iv)] $\mG_2(q)$ for some $q=2^{2e}$, except $\mG_2(4)$.

\end{itemize}
\end{defn}
These are exactly the groups of Lie type in characteristic $2$ that are listed in the table in
Theorem 3 on p. 316 in \cite{KL}.

\begin{thm}
\label{FTKL} Let us assume that $\fchar(K)=0$.
 Suppose that $\End^0(X)$ is a simple $\Q$-algebra with the
center $C$. Suppose that $\tilde{G}_{2,X,K}$ is a simple
non-abelian group, whose order is not a divisor of $2\dim(X)$ and
$\End_{\tilde{G}_{2,X,K}}(X_2)\cong \F_4$. Assume, in addition,
that $\tilde{G}_{2,X,K}$ is a known simple group that is not FTKL-exceptional.
Suppose also that $\dim(X)$
coincides with the smallest positive integers $d$ such that
$\tilde{G}_{2,X,K}$ is isomorphic to a subgroup of $\PGL(d,\C)$.
Then:

\begin{itemize}
\item[(i)] The center $C$ of $\End^0(X)$ is either $\Q$ or
quadratic field.

\item[(ii)] Either  $\End^0(X)=C$ or the following conditions
hold:

\begin{enumerate}
\item There exists a finite algebraic field extension $L/K$ such that
$\tilde{G}_{2,X,L}=\tilde{G}_{2,X,K},$, the overfield
 $L$ is $2$-balanced with respect to $X$, the algebra $E:=\End^0_L(X)$ is either $\Q$
or a quadratic field, $C=E$ and the following conditions hold.

There exist a
finite perfect group $\Pi \subset \End^0(X)^{*}$ and a surjective
homomorphism $\Pi \to \tilde{G}_{2,X,K}$ that is a minimal cover
and a central extension. In addition, the induced homomorphism
$$E[\Pi] \to \End^0(X)$$
is surjective, i.e., $\End^0(X)$ is isomorphic to a direct summand of the group
algebra $E[\Pi]$.

\item  If $C=\Q$ then $X$ enjoys one of the following two
properties:
\begin{itemize}
\item[(a)] $X$ is isogenous over $K_a$ to a self-product of  an
elliptic curve without complex multiplication.
\item[(b)]  $\dim(X)$ is even and $X$ is isogenous over $K_a$ to a
self-product of an abelian surface $Y$ such that $\End^0(Y)$ is an
indefinite quaternion $\Q$-algebra.
\end{itemize}

\item If $C \ne\Q$ then
 $C$ is an imaginary quadratic field and
  $X$ is isogenous over $K_a$ to a
self-product of  an elliptic curve with complex multiplication by
$C$.
\end{enumerate}
\end{itemize}
\end{thm}

\begin{proof}
Using Theorem \ref{center} and replacing if necessary $K$ by its
suitable extension, we may assume that $K$ is $2$-balanced with
respect to $X$, the algebra  $E:=\End^0_K(X)$ is either $\Q$ or a
quadratic field, $C \subset E$ and the following conditions hold.

{\sl Either $\End^0(X)=E$ or there exist a finite perfect group
$\Pi \subset C_{X,K}^{*}$ and a surjective homomorphism $\Pi \to
\tilde{G}_{2,X,K}$ that is a minimal cover and such that the
induced homomorphism
$$E[\Pi] \to C_{X,K}$$
is surjective, i.e., $C_{X,K}$ is isomorphic to a direct summand of the group
algebra} $E[\Pi]$. (Here as above $C_{X,K}$ is the centralizer of $E$ in
$\End^0(X)$). 

Assume that $\End^0(X)\ne E$. We are going to prove that $\Pi \to
\tilde{G}_{2,X,K}$ is a central extension, using results of
Feit-Tits and Kleidman-Liebeck \cite{FT,KL}. Without loss of
generality we may assume that there is a field embedding
$K_a\hookrightarrow \C$ and consider $X$ as complex abelian
variety. Let $t_X$ be the Lie algebra of $X$ that is a
$\dim(X)$-dimensional complex vector space. By functoriality, this
gives us the embeddings 
\begin{eqnarray*}
&&j:\End^0(X)\hookrightarrow \End_{\C}(t_X)\cong \M_{\dim(X)}(\C), \\
&&j:\End^0(X)^*\hookrightarrow \Aut_{\C}(t_X)\cong \GL(\dim(X),\C).
\end{eqnarray*}
Clearly, only central elements of $\Pi$ go to
scalars under $j$. It follow that there exists a central subgroup
$Z$ in $\Pi$ such that $j(Z)$ consists of scalars and
$\Pi/Z\hookrightarrow\PGL(\dim(X),\C)$.
The simplicity of
$\tilde{G}_{2,X,K}$ implies that $Z$ lies in the kernel of $\Pi
\to \tilde{G}_{2,X,K}$ and the induced map $\Pi/Z \to
\tilde{G}_{2,X,K}$ is also a minimal cover. It follows from
Theorem on p. 1092 of \cite{FT} and  Theorem 3 on p. 316 of
\cite{KL} that $\Pi/Z \to \tilde{G}_{2,X,K}$ is a central
extension of $\tilde{G}_{2,X,K}$. Since $\Pi$ is a central
extension of $\Pi/Z $, it follows \cite{Suzuki} that $\Pi$ is a
central extension of $\tilde{G}_{2,X,K}$.

Now notice that $t_X$ carries a natural structure of
$E\otimes_{\Q}\C$-module. Assume that $E\ne\Q$, i.e., $E$ is a
quadratic field. Let $\sigma,\tau: E \hookrightarrow \C$ be the
two different embeddings of $E$ into $\C$. Then
$$E\otimes_{\Q}\C=\C_{\sigma}\oplus\C_{\tau}$$ with
$$\C_{\sigma}=E\otimes_{E,\sigma}\C=\C,\
\C_{\tau}=E\otimes_{E,\tau}\C=\C$$ and
 $t_X$ splits into a direct sum
 $$t_X=\C_{\sigma}t_X\oplus \C_{\tau}t_X.$$
 Suppose that both $\C_{\sigma}t_X$ and $\C_{\tau}t_X$ do not vanish.
 Then the $\C$-dimension $d_{\sigma}$ of non-zero $\C_{\sigma}t_X$ is strictly less
 than dim(X). Clearly, $\C_{\sigma}t_X$ is $C_{X,K}$-stable and we
 get a nontrivial homomorphism
$$C_{X,K}\to \End_{\C}(\C_{\sigma}t_X)\cong \M_{d_{\sigma}}(\C)$$
that must be an embedding in light of the simplicity of $C_{X,K}$.
This gives us an embedding
$$C_{X,K}^{*}\to \Aut_{\C}(\C_{\sigma}t_X)\cong\GL(d_{\sigma},\C).$$
One may easily check that all the elements of $\Pi$ that go to
scalars in $\Aut_{\C}(\C_{\sigma}t_X)$ constitute a central
subgroup $Z_{\sigma}$ that lies in the kernel of $\Pi\to
\tilde{G}_{2,X,K}$. This gives us a central extension
$\Pi/Z_{\sigma}\twoheadrightarrow \tilde{G}_{2,X,K}$ that is a
minimal cover and an embedding $\Pi/Z_{\sigma}\hookrightarrow \Aut_{\C}(\C_{\sigma}t_X)
\cong\GL(d_{\sigma},\C)$. Since $d_{\sigma}<\dim(X)$, Theorem on
p. 1092 of \cite{FT} and  Theorem 3 on p. 316 of \cite{KL} provide
us with a contradiction.
It follows that either $\C_{\sigma}t_X$
or $\C_{\tau}t_X$ does vanish. We may assume that
$\C_{\tau}t_X=0$. This means that each $e\in E$ acts on $t_X$ as
multiplication by complex number $\sigma(e)$, i.e., $j(E)$
consists of scalars. Recall that the exponential map identifies
$X(\C)$ with the complex torus $t_X/\Lambda$ where $\Lambda$ is a
discrete lattice of rank $2\dim(X)$. In addition, $\Lambda$ is
$j(\End_K(X))$-stable where $\End_K(X)$ is an order in the
quadratic field $\End_K^0(X)$.  Now the discreteness of $\Lambda$
implies that $E$ cannot be real and therefore is an imaginary
quadratic field. It follows easily that $X$ is isogenous over $\C$
to a self-product of an elliptic curve with complex multiplication
by $E$. In particular, $E=C$ and $C_{X,K}=\End^0(X)$.

Now let us assume that $E=\Q$. Then $C_{X,K}=\End^0(X)$. Let $Y$
be an absolutely simple abelian variety such that $X$ is isogenous
to a self-product $Y^r$ for some positive integer $r$ with
$r\mid\dim(X)$. Then $\End^0(X)\cong\M_r(\End^0(Y))$.
In
particular, the center of the division algebra $\End^0(Y)$ is
$\Q$. It follows from Albert's classification \cite{MumfordAV}
that $\End^0(Y)$ is either $\Q$ or a quaternion $\Q$-algebra.

If $\End^0(Y)=\Q$ then $\End^0(X)\cong \M_r(\Q)$ and
$\Pi/Z\hookrightarrow \PGL(r,\Q)\subset \PGL(r,\C)$. It follows
that $r=\dim(X)$, i.e. $Y$ is an elliptic curve without complex
multiplication.

Suppose that $\End^0(Y)$ is  a quaternion $\Q$-algebra. Since
$\dim(Y) =\dim(X)/r$ and we live in characteristic zero, $2r$
divides $\dim(X)$. Clearly,
$$\End^0(Y)\subset \End^0(Y)\otimes_{\Q}\C\cong \M_2(\C)$$ and
therefore $$\End^0(X)\cong\M_r(\End^0(Y))\hookrightarrow
\M_{2r}(\C).$$ This implies that $\Pi \hookrightarrow \GL(2r,\C)$.
It follows that $2r=\dim(X)$, i.e., $\dim(Y)=2$. It follows from
the classification of endomorphism algebras of abelian surfaces
\cite[Sect. 6]{Oort} that $\End^0(Y)$ is  an {\sl indefinite}
quaternion $\Q$-algebra.
\end{proof}

\begin{thm}
\label{PSL2} Let us assume that $\fchar(K)=0$.
 Suppose that $\End^0(X)$ is a simple $\Q$-algebra.
 Suppose that $d:=\dim(X)=(q-1)/2$ where $q \ge 5$ is an odd prime power.
Suppose that $\tilde{G}_{2,X,K}\cong \PSL_2(\F_q)$ and
$\End_{\tilde{G}_{2,X,K}}(X_2)\cong \F_4$. Then one of the
following two conditions holds:
\begin{itemize}
 \item[(i)]
$\End^0(X)=\Q$ or a quadratic field. In particular, $X$ is an
absolutely simple abelian variety.

\item[(ii)] $q$ is congruent to $3$ modulo $4$ and $X$ is
$K_a$-isogenous to a self-product of an elliptic curve with
complex multiplication by $\Q(\sqrt{-q})$.
\end{itemize}
\end{thm}

\begin{proof}

 It is well-known \cite[Sect. 4.15]{Gor} that $\SL_2(\F_q)$ is
the universal central extension of $\PSL_2(\F_q)$ and therefore every
projective representation of $\PSL_2(\F_q)$ lifts to a linear representation of
$\SL_2(\F_q)$. The well-known list of irreducible representations of
$\SL_2(\F_q)$ over complex numbers  \cite[Sect. 38]{DA} tells us that the
smallest degree of a nontrivial representation of $\SL_2(\F_q)$ is $(q-1)/2=d$.
This implies that we are in position to apply Theorem \ref{FTKL}. In
particular, $C$ is either $\Q$ or a quadratic field. We may and will assume
that $\End^0(X)\ne C$.

We need to rule out the following possibilities:
\begin{enumerate}
\item\label{case1} $\dim(X)$ is even and $X$ is isogenous over $K_a$ to a
self-product of  an abelian surface $Y$ such that $\End^0(Y)$ is
an indefinite quaternion $\Q$-algebra. In particular, $\End^0(X)$
is a $d^2$-dimensional central simple $\Q$-algebra.

\item\label{case2} $q$ is congruent to $1$ modulo $4$ and
  $X$ is isogenous over $K_a$ to a
self-product of  an elliptic curve with complex multiplication.
In particular, $\End^0(X)$ is a $d^2$-dimensional central simple
algebra over the imaginary quadratic field $C$ unramified at $\infty$.

\item\label{case3} $X$ is isogenous over $K_a$ to a self-product of  an
elliptic curve without complex multiplication. In particular, $\End^0(X)$
is a $d^2$-dimensional central simple $\Q$-algebra.
\end{enumerate}

 By Theorem \ref{FTKL}, there exist a finite perfect
group $\Pi$ and a minimal central cover $\Pi \to \PSL_2(\F_q)$ such that
$\End^0(X)$ is a quotient of the group algebra $E[\Pi]$ where $E=C$ is either
$\Q$ or an imaginary quadratic field. It follows easily that $\Pi=\PSL_2(\F_q)$
or $\SL_2(\F_q)$, so we may always view $\End^0(X)$ as a simple quotient
(direct summand) $D$ of $E[\SL_2(\F_q)]$. By Theorem \ref{FTKL}, $\End^0(X)$ is
a central simple $E$-algebra of dimension $d^2$.

Let us consider the composition
 $$\Q[\SL_2(\F_q)]\subset E[\SL_2(\F_q)] \twoheadrightarrow
 \End^0(X).$$
Let $D$ be the {\sl simple} direct summand of $\Q[\SL_2(\F_q)]$, whose image in
the {\sl simple} $\Q$-algebra $\End^0(X)$ is {\sl not} zero. We write $B\subset
\End^0(X)$ for the image of $D$: it is a $\Q$-subalgebra isomorphic to $D$. The
induced map $D \to \End^0(X)$ is injective, because $D$ is a simple
$\Q$-algebra. On the other hand, $D_E=D\otimes_{\Q}E$ is a direct summand of
$E[\SL_2(\F_q)]$ and the image of $D_E \to \End^0(X)$ is a non-zero ideal of
$\End^0(X)$. Since $\End^0(X)$ is simple, $D_E \to \End^0(X)$ is surjective. In
particular, $B$ generates $\End^0(X)$ as $E$-vector space and the center of $D$
embeds into the center $E$ of $\End^0(X)$. This implies that the center of $D$
is either $\Q$ or isomorphic to $E$. In addition, if the center of $D$ is
isomorphic to $E$ then $B$ contains $E$, i.e., $B$ is a $E$-vector subspace of
$\End^0(X)$ and therefore coincides with $\End^0(X)$: this implies that
$\End^0(X)\cong D$.

Assume that the center of $D$ is isomorphic to $E$. Then  $\End^0(X)\cong D$
and therefore $D$ is a central simple $E$-algebra of dimension $d^2$. This
means that the simple direct summand $D$ of $\Q[\SL_2(\F_q)]$ corresponds to an
irreducible (complex) character of $\SL_2(\F_q)$ of degree $d$ as in Lemma 24.7
of \cite{DA}. These simple direct summands are described explicitly in
\cite{J,F}. In particular, if  $q$ is congruent to $1$ modulo $4$ but is
not a
square then the center of $D$ is a real quadratic field $\Q(\sqrt{q})$, which
is not the case.  This implies that $q$ is congruent to $3$ modulo $4$: in this
case the center of $D$ is an imaginary quadratic field $\Q(\sqrt{-q})$ and
therefore $E=\Q(\sqrt{-q})$. It follows from Theorem \ref{FTKL}  that $X$ is
$K_a$-isogenous to a self-product of an elliptic curve with complex
multiplication by $\Q(\sqrt{-q})$.

Now assume that the center of $D$ is {\sl not} isomorphic to $E$. Then it must
be $\Q$, i.e., $D$ is a central simple $\Q$-algebra. It follows that $D_E$ is a
central simple $E$-algebra and therefore the surjective homomorphism $D_E \to
\End^0(X)$ is injective. It follows that $D_E \cong \End^0(X)$; in particular,
the  central simple $\Q$-algebra $D$ has $\Q$-dimension $d^2$. As above, this
means that the simple direct summand $D$ corresponds to an irreducible
(complex) character of $\SL_2(\F_q)$ of degree $d$. Since the center of $D$ is
$\Q$, it follows from results of \cite{J,F} that $q$ is a square, which is not
the case.
This ends the proof.
\end{proof}

\section{Hyperelliptic jacobians}
\label{mainp}

Suppose that $f(x)\in K[x]$ is a polynomial of degree $n\ge 5$
without multiple roots. Let $\RR_f\subset K_a $ be the set of
roots of $f$. Clearly, $\RR_f$ consists of $n$ elements. Let
$K(\RR_f)\subset K_a$ be the splitting field of $f$. Clearly,
$K(\RR_f)/K$ is a Galois extension and we write $\Gal(f)$ for its
Galois group $\Gal(K(\RR_f)/K)$. By definition, $\Gal(K(\RR_f)/K)$
permutes elements of  $\RR_f$; further we identify $\Gal(f)$ with
the corresponding subgroup of $\Perm(\RR_f)$,
where $\Perm(\RR_f)$
is the group of permutations of $\RR_f$.

We write $\F_2^{\RR_f}$ for the $n$-dimensional $\F_2$-vector
space of maps $h:\RR_f \to \F_2$. The space $\F_2^{\RR_f}$ is
provided with a natural action of $\Perm(\RR_f)$ defined as
follows. Each $s \in \Perm(\RR_f)$ sends a map
 $h:\RR_f\to \F_2$ to  $sh:\alpha \mapsto h(s^{-1}(\alpha))$. The permutation module $\F_2^{\RR_f}$
contains the $\Perm(\RR_f)$-stable hyperplane
$$(\F_2^{\RR_f})^0=
\{h:\RR_f\to \F_2\mid\sum_{\alpha\in \RR_f}h(\alpha)=0\}$$ and the
$\Perm(\RR_f)$-invariant line $\F_2 \cdot 1_{\RR_f}$ where
$1_{\RR_f}$ is the constant function $1$. Clearly,
$(\F_2^{\RR_f})^0$ contains $\F_2 \cdot 1_{\RR_f}$ if and only if
$n$ is even.

 If  $n$ is even then
let us define the $\Gal(f)$-module
$Q_{\RR_f}:=(\F_2^{\RR_f})^0/(\F_2 \cdot 1_{\RR_f})$. If $n$ is
odd then let us put $Q_{\RR_f}:=(\F_2^{\RR_f})^0$. If $n \ne 4$
the natural representation of $\Gal(f)$ is faithful, because in
this case the natural homomorphism
$\Perm(\RR_f)\to\Aut_{\F_2}(Q_{\RR_f})$ is injective.

The canonical surjection $\Gal(K)\twoheadrightarrow
\Gal(K(\RR_f)/K)=\Gal(f)$ provides $Q_{\RR_f}$ with a natural
structure of $\Gal(K)$-module. It is  well-known that the
$\Gal(K)$-modules $J(C_f)_2$ and $Q_{\RR_f}$ are isomorphic (see
for instance \cite{Poonen,SPoonen,ZarhinTexel}). It follows easily
that $K(J(C_f)_2)=K(\RR_f)$ and $\tilde{G}_{2,J(C_f),K}=\Gal(f)$.

Let us put $X=J(C_f)$ and $G := \tilde{G}_{2,X,K}$.
Then $G \cong \Gal(f)$, and the $G$-modules $X_2$ and
$Q_{\RR_f}$ are isomorphic. We freely interchange these two
modules throughout this section.

\begin{ex}
\label{F4c} Suppose that  $n=q+1$ where $q\ge 5$  is a power of an
odd prime $p$. Suppose that $\Gal(f)=\PSL_2(\F_q)$. Assume that
that $f(x)$ is {\sl irreducible}, i.e., $\Gal(f)=\PSL_2(\F_q)$
acts transitively on the $(q+1)$-element set $\RR_f$. If $\beta\in
\RR_f$ then its stabilizer $\Gal(f)_{\beta}$ is a subgroup of
index $q+1$ and therefore contains a Sylow $p$-subgroup of
$\PSL_2(\F_q)$. It follows from the classification of subgroups of
$\PSL_2(\F_q)$ \cite[Theorem 6.25 on page 412]{Suzuki} and
explicit description of its Sylow $p$-subgroup and their
normalizers \cite[p. 191--192]{Huppert} that that
$\Gal(f)_{\beta}$ is conjugate to the (Borel) subgroup of
upper-triangular matrices and therefore the $\PSL_2(\F_q)$-set
$\RR_f$ is isomorphic to the projective line $\P^1(\F_q)$ with the
standard action of $\PSL_2(\F_q)$ (by fractional-linear
transformations),  which is well-known to be doubly transitive \cite{Mortimer}.

Assume, in addition that $q$ is congruent to $\pm 3$ modulo $8$.
Then it is known \cite{Mortimer} that
$$\End_{\Gal(f)}(Q_{\RR_f})=\F_4.$$
\end{ex}

\begin{thm}
\label{l2q} Suppose that $\fchar(K)\ne 2$ and $n=q+1$ where $q\ge
5$ is a prime  power that is congruent to $\pm 3$ modulo $8$.
Suppose that $\Gal(f)=\PSL_2(\F_q)$ acts doubly transitively on
$\RR_f$ (where $\RR_f$ is identified with the projective line
$\P^1(\F_q)$). Then $\End^0(J(C_f))$ is a simple $\Q$-algebra,
i.e. $J(C_f)$ is either absolutely simple or isogenous to a power
of an absolutely simple abelian variety.
\end{thm}

\begin{proof}
See \cite[Theorem 3.10]{ZarhinLuminy}.
\end{proof}

\begin{proof}[Proof of Theorem \ref{main}]
The result follows from Theorem \ref{PSL2} combined with Example
\ref{F4c} and Theorem \ref{l2q}.
\end{proof}

\section{Criteria for Absolute Simplicity}
\label{noniso}

Sometimes, it is possible to rule out the second outcome of
Theorem \ref{main}. First, recall Goursat's lemma \cite[p.
75]{LangAlg}:

\begin{lem}
Let $G_1$ and $G_2$ be finite group, and $H$ a subgroup of $G_1 \times G_2$ such that
the restrictions $p_1: H \to G_1$ and $p_2: H \to G_2$ of the projection maps are surjective.
Let $H_1$ and $H_2$ be the normal subgroups of $G_1$ and $G_2$, respectively,
such the groups $H_1 \times \{1\}$ and $\{1\} \times H_2$ are kernels of $p_2$ and $p_1$,
respectively.
Then there exists an isomorphism $\gamma: G_1/H_1 \cong G_1/H_2$
such that $H$ coincides with the preimage in $G_1\times G_2$ of the
graph of $\gamma$ in $G_1/H_1 \times G_2/H_2$.
\end{lem}

\begin{ex}\label{egprod}
Let $G_1$ be a finite simple group and $G_2$ be a finite group that
does not admit $G_1$ as a quotient. If $H$ is a subgroup of $G_1\times G_2$
that
satisfies the conditions of Goursat's lemma, then $H=G_1 \times G_2$.

Indeed, since $G_1$ is simple, $H_1=\{1\}$ or $G_1$.
We have $H_1\neq\{1\}$, since otherwise $G_1/H_1 \cong G_1$ and
no quotient of $G_2$ is isomorphic to $G_1$.
Therefore, $H_1 = G_1$, $G_1/H_1 \cong G_2/H_2 = \{1\}$,
and $H_2=G_2$.
Since $G_1/H_1\times G_2/H_2$ is a trivial group,
the graph of $\gamma$ coincides with $G_1/H_1\times G_2/H_2$,
and its preimage $H$ coincides with $G_1 \times G_2$.
\end{ex}

\begin{thm}\label{thmnoniso}
Let $K$ be a field of characteristic zero. Suppose that  $f(x)\in
K[x]$ is a  polynomial of degree $n \ge 5$ without multiple roots.
Let us consider the hyperelliptic curve $C_f: y^2 = f(x)$ and its
jacobian $J(C_f)$. Suppose that $h(x) \in K[x]$ is an irreducible
cubic polynomial and let us consider the elliptic curve $Y: y^2 =
h(x)$.

Let us  assume that $f(x)$ and $h(x)$ enjoy the following
properties:

\begin{enumerate}
\item\label{simplecond}
$\Gal(K(\RR_f)/K) = \PSL_2(\F_q)$ for some odd prime power $q\equiv 3 \mod 8$
with $n=q+1$, and $\Gal(K(\RR_f)/K)$ acts doubly transitively on
$\RR_f$ (where $\RR_f$ is identified with the projective line $\P^1(\F_q)$);
\item\label{s3cond} $\Gal(K(\RR_h)/K) = \SS_3$.
\end{enumerate}
Then $\Hom(J(C_f), Y)=0$ and $\Hom(Y, J(C_f))=0$. In particular,
$J(C_f)$ is not $K_a$-isogenous to a self-product of $Y$.
\end{thm}
\begin{proof}
First, we prove that $K(\RR_f)$ and $K(\RR_h)$ are linearly
disjoint over $K$. Let us put $G_1:=\Gal(K(\RR_f)/K)$,
$G_2:=\Gal(K(\RR_h)/K)$, and $H := \Gal(K(\RR_f, \RR_h)/K)$, the
Galois group of the compositum of $K(\RR_f)$ and $K(\RR_h)$ over
$K$. By Theorem 1.14 of \cite{LangAlg}, $H$ can be considered to
be a subgroup of $G_1\times G_2$, where the Galois restriction
maps coincide with restrictions of projection maps $p_i: G_1
\times G_2 \to G_i$, with $i=1, 2$, to $H$. It follows from
Example \ref{egprod} that $H \cong G_1 \times G_2$, and $K(\RR_f)$
and $K(\RR_h)$ are linearly disjoint over $K$. The equalities
 $\Hom(J(C_f), Y)=0$ and $\Hom(Y, J(C_f))=0$
 follow from the definitions (s) and (p3) and Theorem 2.5 of
\cite{ZarhinHHJ}. Since for any positive integer $r$ we have
$\Hom(J(C_f), Y^r)= \prod_{i=1}^r \Hom(J(C_f),Y)$,
we conclude that $\Hom(J(C_f), Y^r) = 0$.
\end{proof}

The following assertion will be proven in Section \ref{quadratic}.

\begin{thm}
\label{class} Let $p>3$ be a prime such that $p\equiv 3 \mod 8$.
Let us put $\omega=\frac{-1+\sqrt{-p}}{2}$ and let
$\O=\Z+\Z\omega$ be the ring of integers in $\Q(\sqrt{-p})$. Let
$\O_2=\Z+2\O$ be the order of conductor $2$ in $\Q(\sqrt{-p})$.
\begin{itemize}
\item[(i)] The principal ideal $(2)$ is prime in $\O$.
 \item[(ii)] Let $\b$ be a proper fractional $\O_2$-ideal in
$\Q(\sqrt{-p})$ and $\a=\O\b$ be
 the $\O$-ideal generated by $\b$.
Then $\b$ contains $2\a$ as a subgroup of index $2$ and $\a$
contains $\b$ as a subgroup of index $2$.
\item[(iii)] Let $\a$ is
a  fractional $\O$-ideal in $\Q(\sqrt{-p})$. If $\b$ is a subgroup
of index $2$ in $\a$ then it is a proper $\O_2$-ideal in
$\Q(\sqrt{-p})$, i.e.,
$$\O_2=\{z\in \Q(\sqrt{-p})\mid z\b \subset \b\};$$
in addition, $\a=\O\b$.
 There are exactly three index $2$ subgroups in $\a$;
they are mutually non-somorphic as $\O_2$-ideals.

\item[(iv)] If $\hh$ is the class number of $\Q(\sqrt{-p})$ then
$3\hh$ is the number of classes of proper $\O_2$-ideals.
\end{itemize}
\end{thm}

We write $j$ for the classical modular function  \cite[Ch. 3,
Sect. 3]{LangE}.

\begin{cor} \label{jinvcor}
Let $p$ be a prime such that $p\equiv 3 \mod 8$.
Let $q\ge 11$ be an odd power of $p$. (In particular, $q
\equiv p\equiv 3 \mod 8$.) Let us put
$$\omega:=\frac{-1+\sqrt{-p}}{2}, \ \alpha: = j(\omega)\in \C, \ K:=\Q(j(\omega))\subset \C.$$
Suppose that $f(x)\in K[x]$ is an irreducible polynomial of degree
$q+1$ such that $\Gal(f/K)=\PSL_2(\F_q)$ acts doubly transitively
on $\RR_f$ (where $\RR_f$ is identified with the projective line
$\P^1(\F_q)$).

Then $J(C_f)$ is an absolutely simple abelian variety, and
$\End^0(J(C_f))=\Q$ or a quadratic field.
\end{cor}

\begin{proof}
Clearly, $\Q(\sqrt{-p})=\Q(\sqrt{-q})$. Since $p\equiv 3 \mod 8$,
the ring of integers $\O$ in  $\Q(\sqrt{-p})$ coincides with
$\Z+\Z\omega$. If $p>3$ let us consider the polynomial
$$
h_p(x):=x^3 - \dfrac{27\alpha}{4(\alpha-1728)}x -
\dfrac{27\alpha}{4(\alpha-1728)}\in K[x].
$$
If $p=3$ then $$\omega=\frac{-1+\sqrt{-3}}{2},\ \alpha =
j(\omega)=0,\ K=\Q$$ and we put
$$h_3(x):=x^3-2\in \Q[x]=K[x].$$

The elliptic curve $Y: y^2=h_p(x)$ is defined over $K$ and its
$j$-invariant coincides with $\alpha=j(\omega)$, i.e.,
$Y(\C)=\C/(\Z+\Z\omega)$. Hence $Y$ admits complex multiplication
by $\Z+\Z\omega=\O$.

The following Lemma will be proven at the end of this Section.
\begin{lem}
\label{GalS3}
 The polynomial $h_p(x)$ is irreducible over $K$ and its
Galois group $\Gal(h_p/K) \cong \SS_3$.
\end{lem}

Combining Lemma \ref{GalS3} and Theorem \ref{thmnoniso}, we
conclude that $J(C_f)$ is not isogenous to a self-product of $Y$.
Now the result follows from Theorem \ref{main}.
\end{proof}

\begin{proof}[Proof of Lemma \ref{GalS3}]
The case $p=3$ is easy. So, further we assume that $p>3$. First,
check that  the discriminant $\Delta$ of $h_p$ is not a square in
$K$. Indeed, we have $\Delta = 1458^2\, \alpha^2/(\alpha-1728)^3$,
so $\Delta$ is a square in $K$ if and only if $\alpha-1728$ is.
According to \cite[p. 288]{Birch}, if $p\equiv 3 \mod 4$, then
$\alpha$ is real and negative. Therefore the purely imaginary
$\sqrt{\alpha-1728}$ does not lie in the real number field $K :=
\Q(\alpha)$. Now it suffices to check that the cubic polynomial
$h_p$ is irreducible. Suppose that this is not the case, i.e.,
$h_p$ has a root in $K$. This means that $Y$ has a $K$-rational
point of order $2$ and therefore there exists an elliptic curve
$Y'$ over $K$ and a degree $2$ isogeny $Y \to Y'$. By duality, we
get a degree $2$ isogeny $Y' \to Y$. This allows us to identify
$Y'(\C)$ with $\C/\b$ where $\b$ is a subgroup of index
 $2$ in $\Z+\Z\omega=\O$. By Theorem \ref{class}(ii), $\b$ is a proper
 $\O_2$-ideal. The classical theory of complex multiplication
 \cite[Th. 5.7 on p. 123]{Shimura} tells us that
 $\Q(\sqrt{-p})(j(Y'))/\Q(\sqrt{-p})$ is an abelian extension,
 whose degree coincides with the  number of
classes of proper $\O_2$-ideals. By Theorem \ref{class}(iv),
$[\Q(\sqrt{-p})(j(Y')):\Q(\sqrt{-p})]=3\hh$. On the other hand,
since $Y'$ is defined over $K$, the number $j(Y')\in K$ and
therefore
$$\Q(\sqrt{-p})(j(Y'))\subset
\Q(\sqrt{-p})K=\Q(\sqrt{-p})(\alpha)=\Q(\sqrt{-p})(j(Y)).$$
However, $\Q(\sqrt{-p})(j(Y))$ is the absolute class field of
$\Q(\sqrt{-p})$ and it is well-known \cite{SerreCM,Birch},
\cite[Th. 5.7 on p. 123]{Shimura} that
$[\Q(\sqrt{-p})(j(Y)):\Q(\sqrt{-p})]=\hh$. Since $\hh<3\hh$, we
obtain the desired contradiction.
\end{proof}


\begin{rem}
If $p=3,11, 19, 43, 67$ or $163$ then  $j(\omega)\in \Z$
\cite{SerreCM}, so  $K =\Q(\alpha)= \Q(j(\omega))=\Q$.
\end{rem}

\begin{thm} \label{Qjinvcor}
Let $p$ be a prime such that $p\equiv 3 \mod 8$.
Let $q\ge 11$ be an odd power of $p$. (In particular, $q \equiv
p\equiv 3 \mod 8$.) Suppose that $f(x)\in \Q[x]$ is an irreducible
polynomial of degree $q+1$ such that $\Gal(f/\Q)=\PSL_2(\F_q)$
acts doubly transitively on $\RR_f$ (where $\RR_f$ is identified
with the projective line $\P^1(\F_q)$).

Then $J(C_f)$ is an absolutely simple abelian variety, and
$\End^0(J(C_f))=\Q$ or a quadratic field.
\end{thm}

\begin{proof}
Let us put
$$\omega:=\frac{-1+\sqrt{-p}}{2}, \ \alpha: = j(\omega)\in \C, \ K:=\Q(j(\omega))\subset \C.$$
Since simple non-abelian $\PSL_2(\F_q)$ does not have a subgroup
of index $2$,
 $$\Gal(f/\Q)=\PSL_2(\F_q)=\Gal(f/\Q(\sqrt{-p})).$$
 Since $\PSL_2(\F_q)$ is perfect and
 $K\Q(\sqrt{-p})=\Q(\sqrt{-p})(j(\omega))$ is abelian over
 $\Q(\sqrt{-p})$,
$$\Gal(f/\Q)=\PSL_2(\F_q)=\Gal(f/K\Q(\sqrt{-p})).$$
Since $$\Gal(f/K\Q(\sqrt{-p}))\subset \Gal(f/K)\subset
\Gal(f/\Q),$$ we conclude that
$$\Gal(f/K)=\Gal(f/\Q)=\PSL_2(\F_q)$$
acts doubly transitively on $\RR_f$. In order to finish the proof,
one has only to apply Corollary \ref{jinvcor}.
\end{proof}

\section{Proof of Theorem \ref{class}}
\label{quadratic} There is a positive integer $k$ such that
$p=8k+3$. It follows that $\omega^2+\omega+(2k+1)=0$. This implies
that the $4$-element algebra $\O/2\O$ contains a subalgebra
isomorphic to the finite field $\F_4$ and therefore coincides with
$\F_4$. This means that $(2)$ is prime in $\O$. So, this proves
(i).

Suppose that $\b$ is a proper $\O_2$-ideal in $\Q(\sqrt{-p})$ and
$\a:=\O\b$. Clearly, $2\a\subset\b\subset\a$. Since $\a$ and $2\a$
are $\O$-ideals, $\b$ does coincides neither with $\a$ nor with
$2\a$. Since $2\a$ has index $4$ in $\a$, the group $\b$ has index
$2$ in $\a$ and $2\a$ has index $2$ in $\b$. This proves (ii).

Now, suppose that $\a$ is a fractional $\O$-ideal in
$\Q(\sqrt{-p})$ and a subgroup $\b\subset \Q(\sqrt{-p})$ satisfies
$2\a\subset\b\subset\a$. If $\b$ is an $\O$-ideal then the unique
factorization of $\O$-ideals and the fact that $(2)$ is prime
imply that either $\b=\a$ or $\b=2\a$. So, if $\b$ has index $2$
in $\a$, it is neither $\a$ nor $2\a$ and therefore is {\sl not}
an $\O$-ideal.

On the other hand, it is clear that $\O\b\subset\a$ and
$2\O\b\subset 2\a\subset\b$ and therefore $\b$ is a proper
$\O_2$-ideal. This proves the first assertion of (iii). We have
$\b\subset\O\b\subset\O\a=\a$ but $\b\ne\O\b$. Since the index of
$\b$ in $\a$ is $2$, we conclude that $\O\b=\a$. This proves the
second assertion of (iii).

 Since $\a$ is a free commutative
group of rank $2$, it contains exactly three subgroups of index
$2$. Let $\b_1$ and $\b_2$ be two distinct subgroups of index $2$
in $\a$. We have
 $$\O\b_1=\a=\O\b_2.$$
Suppose that $\b_1$ and $\b_2$ are isomorphic as $\O_2$-ideals.
This means that there exists a non-zero $\lambda\in \Q(\sqrt{-p})$
such that $\lambda \b_1=\b_2$. It follows that $\lambda\a=\a$ and
therefore $\lambda$ is a unit in $\O$. Since $p>3$, we have
$\lambda=\pm 1$ and therefore $\b_2=\b_1$. This proves the last
assertion of (iii).

The assertion (iv) follows easily from (ii) and (iii). (It is also
a special case of Exercise 11 in Sect. 7 of Ch. II in \cite{BS}
and of Exercise 4.12  in Section 4.4 of \cite{Shimura}).


\section{Examples}
\label{examples}

\begin{ex}\label{ex11}
Let $S$ be a transcendental over $\Q$,  $T=2^8 3^5/(11S^2+1)$, and put
\begin{align*}
f_{11, S}(x) :=
& (x^3 - 66x - 308)^4\\
& - 9  T  (11x^5 - 44x^4 - 1573x^3 + 1892x^2 + 57358x + 103763) \\
& - 3  T^2  (x - 11).
\end{align*}
According to Table 10 of the Appendix in \cite{MMIGT},
$\Gal(f_{11, S}/\Q(S)) = \PSL_2(\F_{11})$. It can be verified
using MAGMA \cite{MAGMA} that when $s = m / n$ for any nonzero
integers $-5\leq m,n \leq 5$, then $\Gal(f_{11, s}/\Q) =
\PSL_2(\F_{11})$. Consider the hyperelliptic curve
$$
C_{11, s}: y^2 = f_{11, s}(x)
$$
over $\Q$ by any one of these $s$.
By Theorem \ref{Qjinvcor} the $5$-dimensional abelian variety
$J(C_{11, s})$ is absolutely simple. For example, if we put $s=1$,
then we obtain a hyperelliptic curve
\begin{eqnarray*}
C_{11, 1}: y^2 &=& x^{12} - 264  \,x^{10} - 1232  \,x^9 + 26136  \,x^8 + 243936  \,x^7 - 580800  \,x^6 \\
& &\ \ \ - 16612992  \,x^5 - 54104688  \,x^4 + 310712512  \,x^3 + 2391092352  \,x^2 \\
& &\ \ \ + 4956865152  \,x + 5044849216
\end{eqnarray*}
over $\Q$  with $J(C_{11, 1})$ absolutely simple.
\end{ex}

\begin{ex}\label{ex13}
If we define
\begin{align*}
f_{13, S}(x) :=
& (x^2 + 36)  (x^3 - x^2 + 35x - 27)^4 \\
& - 4  T  (7x^2 - 2x + 247)  (x^2+39)^6 / 27
\end{align*}
with $T=1/(39S^2+1)$ then again \cite{MMIGT} we have
$\Gal(f_{13, S}/\Q(S)) = \PSL_2(\F_{13})$.
Similarly, we checked using MAGMA
 that when $s = m / n$
for any nonzero integers $-5\leq m,n \leq 5$,
then $\Gal(f_{13, s}/\Q) = \PSL_2(\F_{13})$.
If we define
$$
C_{13, s}: y^2 = f_{13, s}(x)
$$
over $\Q$, then by Theorem \ref{FTKL} the $6$-dimensional abelian
variety $J(C_{13, s})$ is absolutely simple. As an example, take $s=-1$ to get
the hyperelliptic curve
\begin{eqnarray*}
C_{13, -1}: y^2 &=& 263/270  \,x^{14} - 539/135  \,x^{13} + 9451/54  \,x^{12} - 10114/15  \,x^{11} \\
&& \ \ \ + 376363/30  \,x^{10}  - 45487 \,x^9 + 891605/2 \,x^8 - 1533844 \,x^7 \\
&& \ \ \ + 15279043/2 \,x^6 - 25943931 \,x^5 + 391472991/10 \,x^4 \\
&& \ \ \ - 896502438/5 \,x^3 - 780396201/2 \,x^2 - 365687757/5 \,x \\
&& \ \ \ - 31998670461/10
\end{eqnarray*}
defined over $\Q$ with $J(C_{13, -1})$ absolutely simple.
\end{ex}

See \cite{Malle} for other examples of irreducible polynomials
over $\Q(T)$ of degrees $n=p+1$ with $p=11, 13, 19, 29,37$, whose
Galois groups are isomorphic to $\PSL_2(\F_p)$. These polynomials
can be used in a manner similar to that of Examples \ref{ex11} and
\ref{ex13}, in order to construct examples of absolutely simple
abelian varieties over $\Q$ of dimensions $5, 6, 9, 14,18$
respectively, whose endomorphism algebra is either $\Q$ or a
quadratic field.


\end{document}